\documentclass[11pt]{amsart}
\usepackage[english]{babel}
\usepackage{amssymb, amsmath, dsfont, epsfig, psfrag}

\newcommand{\R}{\mathds{R}}
\newcommand{\C}{\mathds{C}}

\newcommand{\T}{\mathds{T}}

\renewcommand{\geq}{\geqslant}
\newcommand{\N}{\mathds{N}}

\newcommand{\Z}{\mathds{Z}}

\newcommand{\PP}{\mathds{P}} 
\newcommand{\EE}{\mathds{E}}
\renewcommand{\le}{\leqslant}

\newtheorem{theo}{Theorem}[section]
\newtheorem{ineq}[theo]{Inequality}%[section]
\newtheorem{cor}[theo]{Corollary}%[section]

\newtheorem{rem}{Remark}[section]
\newtheorem{ex}{Example}[section]
\long\def\symbolfootnote[#1]#2{\begingroup\def\thefootnote{\fnsymbol{footnote}}\footnote[#1]{#2}\endgroup}
\newenvironment{dem}
{\if\par\fi
  \rmfamily \noindent
  {\bf Proof\/}~:}{\hfill $\Box$ \par \noindent}

\title
{Random series of functions and applications}

\author {Fr\'ed\'eric Paccaut$^{(1)}$, Dominique Schneider$^{(2)}$}

\begin{document}

\maketitle
{\tiny \noindent(1) Universit\'e de Picardie Jules Verne, L.A.M.F.A. CNRS UMR 6140\\
33, rue Saint Leu, F-80039 Amiens cedex 01,
\texttt{\textbf{frederic.paccaut@u-picardie.fr}}\\
(2) Universit\'e du Littoral C\^ote d'Opale, L.M.P.A. CNRS EA 2597\\
50, rue F.Buisson B.P. 699, F-62228 Calais cedex, 
\texttt{\textbf{dominique.schneider@lmpa.univ-littoral.fr}}}

\vskip3cm
\begin{abstract}
In this article we study the continuity properties of trajectories for some random series of
functions, $\sum_{k=0}^{\infty} a_k f(\alpha X_k(\omega))$ where 
$(a_k)_{k\geq 0}$ is a complex sequence, $(X_k)_{k\geq 0}$ is a sequence of 
real independent random variables, $f$ is a real valued function with period one and summable Fourier 
coefficients. We obtain almost sure continuity results for these periodic or almost periodic series for a large class of functions $f$,
where the "almost sure" does not depend on the function.
The proof relies on gaussian randomization. We show optimality of the results in some cases.
\end{abstract}
\vskip2cm

\begin{center}
{\bf\normalsize SERIES DE FONCTIONS ALEATOIRES ET APPLICATIONS}
\end{center}
\vskip3cm
\begin{abstract}(\textsc{R\'esum\'e})
Dans ce travail, nous \'etudions des propri\'et\'es de continuit\'e de trajectoires de s\'eries de fonctions al\'eatoires du type
$\sum_{k=0}^{\infty} a_k f(\alpha X_k(\omega))$ o\`u $(a_k)_{k\geq 0}$ est une suite de nombres complexes, $(X_k)_{k\geq 0}$ une suite de
variables al\'eatoires r\'eelles et ind\'ependantes, $f$ une fonction 1-p\'eriodique \`a coefficients de Fourier sommables. Nous montrons
que, presque s\^urement, ces s\'eries de fonctions al\'eatoires (p\'eriodiques ou presque p\'eriodiques) sont \`a trajectoires
continues pour une grande classe de fonctions $f$.
Le "presque s\^ur" est ind\'ependant de $f$. Les preuves s'appuient sur un proc\'ed\'e de randomisation gaussien. Dans certains cas, nous
montrerons l'optimalit\'e des r\'esultats obtenus.
\end{abstract}

\symbolfootnote[0]{\hskip-0.4cm{\bf Keywords} : random fourier series, almost sure continuity of trajectories, gaussian randomization,
almost periodic functions, random trigonometric polynomial}

\newpage

\section{Introduction. Main results}
 
In \cite{Berkes}, Berkes studies the almost sure convergence of series defined by :
 $$\sum_{k\geq 1} a_kf(\alpha\, n_k)$$
where the sequence $(n_k)$ is lacunary and the function $f$ verifies :
  $$f(x+1)=f(x)\ \ \ \ \ \int_0^1f(x)dx=0\ \ \ \ \ \int_0^1f^2(x)dx=1$$
He shows that the important property of $f$ to ensure the almost sure convergence is
$f\in\mbox{Lip}(\gamma)$ with $\gamma > 1/2$ and
$\sum_{k\geq 1}|a_k|^2 < +\infty$. In his case, the $n_k$ are strictly lacunary, more precisely, they satisfy
the Hadamard gap condition :
  $$\frac{n_{k+1}}{n_k}\geq q>1$$
We can naturally adress the question whether the convergence still holds when $n_k$ is polynomial, and for which class of functions.
We are going to answer the question when the sequence $(n_k)$ is randomly generated.\\
Let us mention that the result exists when $(n_k)$ is a deterministic polynomial sequence and $(a_k)$ is randomly distributed 
(see \cite{Kahane} and \cite{Pisier}).\\ 
We want to study the convergence properties of series of functions sampled by a random process.
More precisely, consider the torus
${\T} = {\R} / {\Z}$ and define $A({\T})$ as the set of complex valued functions
whose Fourier coefficients are absolutely summable : 
  $$A({\T})=\{f:{\T}\to{\C},f(\alpha)
=\sum_{j\in {\Z}} \hat {f}(j) \exp {(2i\pi \alpha j)},\sum_{j\in {\Z}}
\vert \hat {f}(j) \vert < + \infty\}$$
$(a_k)_{k\geq 0}$ will denote a sequence of real numbers and $(X_k)_{k\geq 0}$ a sequence of
independent real random variables defined on the probabilised space $(\Omega, \mathcal A, \mathbb P).$ \\
Our aim is to study the convergence, when $\omega \in \Omega$ is fixed, of the series of functions
  $$F(\alpha, \omega) =\sum_{k=0}^{\infty} a_k f(\alpha X_k(\omega))\,\leqno {\forall\alpha\in{\R}, }$$
Is it possible to give conditions on the sequence $( a_k)_{k\geq 0}$ in order to find a $\mathcal A-$measurable
set $\Omega_0$ independent of the function $f$, such that $\mathbb P(\Omega_0)=1$,
on which the series uniformly converges?\\
Note that when $X_k$ does not take integer values, $F$ is not a periodical function of the torus.
For us, $\alpha$ will be real and we will deal with this "almost periodical" case. That is why we have to study the properties
of $F$ on a compact $[-M,M]$ and not only $[0,1]$ (see for example \cite{FanSchnei})\\
For all $f\in A({\T})$, define 
  $$\vert \vert  f\vert \vert := \sum_{j\in {\Z}} \vert \hat {f} (j) \vert <+\infty  \, .$$
  $$\vert \vert \vert f \vert \vert \vert := \sum_{j\in {\Z}} \vert \hat {f} (j)
  \vert \sqrt {\log {(\vert j \vert +3)}} < +\infty  \, .$$
and 
  $$B({\T})=\{f:{\T}\to{\C},f(\alpha)
  =\sum_{j\in {\Z}} \hat {f}(j) \exp {(2i\pi \alpha j)}, \vert \vert \vert f \vert \vert \vert < + \infty\}$$
Remark that $B({\T})\subset A({\T})$.\\ 
In the following, we will give conditions for $\alpha\mapsto F(\alpha, \omega)$ to have continuous 
trajectories $\mathbb P-$almost surely.\\
We will denote by $\varphi_{X}$ the characteristic function of the random variable $X$
  $$\forall t\in{\R}, \varphi_{X}(t)={\EE}(e^{2i\pi tX})$$

\begin{theo}\label{moment}
Let $(X_k)_{k\geq 0}$ be a sequence of independent real valued random variables and let 
$(a_k)_{k\geq 1}$ be a sequence of complex numbers such that, for any compact $K$ which does not contain $0$: 
$$
\forall\varepsilon>0, \exists N>0, \sup_{m>n\geq N}\sup_{\alpha\in K}\sup_{j\in {\Z} -\{0\}}\left\vert 
\sum_{k=n}^m a_k\varphi_{X_k}(j\alpha)\right\vert<\varepsilon\, .
\leqno {(\mathcal H)}
$$
Assume moreover that:\\
{\bf case 1:(polynomial)} there exists $\beta>0$ and $d>0$ with ${\EE} \vert X_k \vert^{\beta} = {\mathcal O}(k^d)$ and
\begin{equation}\label{condition1}
  \sum_{n\geq 1}\frac{\sqrt{\sum_{k\geq n}\vert a_k\vert^2}}{n\sqrt{\log n}}<+\infty
\end{equation}
{\bf case 2:(subexponential)} there exists $\beta>0$ and $\gamma\in]0,1[$ with ${\EE} \vert X_k \vert^{\beta} = {\mathcal O}(2^{k^{\gamma}})$ and
\begin{equation}\label{condition2}
  \sum_{n\geq 1}\frac{\sqrt{\sum_{k\geq n}\vert a_k\vert^2}}{n^{1-\frac{\gamma}{2}}}<+\infty
\end{equation}
then in both cases, there exists a measurable set $\Omega_0$ with $\mathbb P(\Omega_0)=1$
such that for all $\omega \in
\Omega_0$, for any $f\in B({\T})$ such that $\int_{{\T}}f(t)dt=0$ : 
for $\alpha\in {\R}-\{0\}$, $F(\alpha, \omega)$ is well defined, $\alpha
\mapsto F(\alpha, \omega)$ is continuous and the series
defining $F$ converges uniformly on every compact which does not contain $\{0\}$.
\end{theo}

\begin{rem}\ \\
\begin{enumerate}
\vspace{-0.5cm}
\item It is worth noticing that the set $\Omega_0$ does not depend on the class of functions $f$ ($A({\T})$ or $B({\T})$).
\item when $(X_k)_{k\geq 0}$ takes integer values, condition
$\vert\vert\vert f\vert\vert\vert < \infty $ becomes $\vert\vert f\vert\vert <\infty $.
\item we will give conditions on the law of the process $(X_k)_{k\geq 0}$ to fulfill hypothesis 
$(\mathcal H)$.
\item For example, when $\vert a_k \vert ={\mathcal O}(k^{-\delta})$,
in case 1, if $\delta>1/2$, then condition \ref{condition1} holds
and in case 2, if $\delta>\frac{\gamma+1}{2}$, then condition \ref{condition2} holds.
\item Concerning case 2, if $\gamma\geq 1$ (${\EE} \vert X_k \vert^{\beta}$ growths exponentially), one can prove using remark \ref{remarque}
that the series $\sum a_k$ has to converge. The function $F$ is then obviously well defined using only Cauchy Schwarz inequality.
\end{enumerate}
\end{rem}

In case condition ($\mathcal H$) is hard to check, it is possible to split up the hypothesis on the sequence $(a_k)$ and the characteristic
function $\varphi_{X_k}$ either using  Abel's summation method or using Cauchy Schwarz inequality.\\
Define: 
  $$c_n=\left\{\begin{array}{cl}
    1+\sqrt{\log n} & \mbox{in the polynomial case} \\
    n^{\frac{\gamma}{2}} & \mbox{in the subexponential case}
    \end{array}\right.$$

\begin{cor}\label{separation}
Let $(X_k)_{k\geq 0}$ be a sequence of independent real valued random variables
\par
Assume that, for any compact $K$ which does not contain $0$ : 
$$
\sup_{N\geq 1}\sup_{\alpha\in K}\sup_{j\in {\Z} -\{0\}}\vert 
\sum_{k=0}^N \varphi_{X_k}(j\alpha)\vert<\infty\, .
\leqno {(\mathcal H')}
$$
Let $(a_k)_{k\geq 1}$ be a sequence of complex numbers enjoying the following properties
\begin{enumerate}
\item $\sum_{n\geq 1}\frac{\sqrt{\sum_{k\geq n}\vert a_k\vert^2}}{nc_n}<+\infty$
\item $\sum_{k\geq 1}\vert a_k-a_{k+1}\vert$ converges
\end{enumerate}
then there exists a measurable set $\Omega_0$ with $\mathbb P(\Omega_0)=1$
such that for all $\omega \in
\Omega_0$, for any $f\in B({\T})$ such that $\int_{{\T}}f(t)dt=0$: 
for $\alpha\in {\R}-\{0\}$, $F(\alpha, \omega)$ is well defined, $\alpha
\mapsto F(\alpha, \omega)$ is continuous and the series
defining $F$ converges uniformly on every compact which does not contain $\{0\}$.
\end{cor}

\begin{cor}\label{separationcauchy}
Let $(X_k)_{k\geq 0}$ be a sequence of independent real valued random variables
\par
Assume that, for any compact $K$ which does not contain $0$ : 
$$
\forall\varepsilon>0, \exists N>0, \sup_{m>n\geq N}\sup_{\alpha\in K}\sup_{j\in {\Z} -\{0\}}\left(
\sum_{k=n}^m \vert\varphi_{X_k}(j\alpha)\vert^2\right)<\varepsilon\, .
\leqno {(\mathcal H'')}
$$
Let $(a_k)_{k\geq 1}$ be a sequence of complex numbers enjoying
  $$\sum_{n\geq 1}\frac{\sqrt{\sum_{k\geq n}\vert a_k\vert^2}}{nc_n}<+\infty$$
then there exists a measurable set $\Omega_0$ with $\mathbb P(\Omega_0)=1$
such that for all $\omega \in
\Omega_0$, for any $f\in B({\T})$ such that $\int_{{\T}}f(t)dt=0$: 
for $\alpha\in {\R}-\{0\}$, $F(\alpha, \omega)$ is well defined, $\alpha
\mapsto F(\alpha, \omega)$ is continuous and the series
defining $F$ converges uniformly on every compact which does not contain $\{0\}$.
\end{cor}  

\begin{rem}\ \\
The previous corollaries will be useful for example when the law of $X_k$ is obtained by convolution product (see corollary \ref{aperiodic}).
As the condition $\sum_{k\geq 1}\vert a_k-a_{k+1}\vert$ is often hard to check, corollary \ref{separationcauchy} is sometimes better to use.
\end{rem}

The proof of theorem \ref{moment} will start by looking separately at $F(\alpha,\omega)-{\EE}(F(\alpha,.))$ and ${\EE}(F(\alpha,.))$. It turns out that
hypothesis $({\mathcal H})$ will be used only to deal with the expectation. That is why we think interesting to state the result for: 
  $$F(\alpha,\omega)-{\EE}(F(\alpha,.)):=\sum_k a_k[f(\alpha X_k(\omega)-{\EE}(f(\alpha X_k))$$ 

\begin{theo}\label{centered}
Let $(X_k)_{k\geq 0}$ be a sequence of independent real valued random variables
such that there exists $\beta>0$ and $d>0$ with ${\EE} \vert X_k \vert^{\beta} = {\mathcal O}(k^d)$ or $\gamma\in]0,1[$ with 
${\EE} \vert X_k \vert^{\beta} = {\mathcal O}(2^{k^{\gamma}})$.
Let $(a_k)_{k\geq 1}$ be a sequence of complex numbers  enjoying the following property
 $$\sum_{n\geq 1}\frac{\sqrt{\sum_{k\geq n}\vert a_k\vert^2}}{nc_n}<+\infty$$
then there exists a measurable set $\Omega_0$ with $\mathbb P(\Omega_0)=1$
such that for all $\omega \in
\Omega_0$, for any $f\in B({\T})$ such that $\int_{{\T}}f(t)dt=0$ : 
for $\alpha\in {\R}$, $F(\alpha,\omega)-{\EE}(F(\alpha,.))$ is well defined, $\alpha
\mapsto F(\alpha,\omega)-{\EE}(F(\alpha,.))$ is continuous, the series
defining $F-{\EE}(F)$ converges uniformly on every compact and there exists $C_{\omega}>0$ such that for all $\alpha\in {\R}$:
  $$\vert F(\alpha,\omega)-{\EE}(F(\alpha,.))\vert\le C_{\omega}\vert\vert\vert f\vert\vert\vert\sqrt{\log(\vert \alpha\vert+2)}$$
\end{theo}

\begin{rem}\ \\
\begin{enumerate}
\vspace{-0.5cm}
\item We also discuss the optimality of hypothesis on $(a_k)$ of theorem \ref{centered} in section 2.
\item We also have:
  $${\EE}\sup_{T>1}\frac{\sqrt{\int_0^T\vert F(t,\omega)-{\EE}(F(t,.))\vert^2dt}}{\sqrt{T\log T}}<\infty$$
\end{enumerate}
\end{rem}

This result relies on uniform estimations of the size of some trigonometric polynomials, 
more precisely on the following :\\
Recall that $\log^+=\max(\log, 0)$.

\begin{theo}\label{principal}
Let $\lambda$ and $\Lambda$ be two integers with $\lambda\le\Lambda$, 
$(X_k)_{k\geq 0}$ be a sequence of independent real valued random variables such that there exists $\beta>0$ such that, $\forall N\geq 0$, 
${\EE}\vert X_N^{\beta}\vert<\infty$. Define
  $$\forall N\geq 0, \Phi_{\beta}(N)=2+\max(N,{\EE}\vert X_N^{\beta}\vert)$$
Let $M\geq 1$ and $I_M = [-M, M].$ Let
$(a_k)_{k\geq 1}$ be a sequence of real or complex numbers.
 \par
Define
  $$A_{\lambda,\Lambda,M}=\sqrt{\log{(M\Phi_{\beta}(\Lambda))}\sum_{k=\lambda}^{\Lambda}\vert a_k\vert ^2},$$
then
  $${\EE}\sup_{j\in {\Z}}\sup_{\lambda\geq 1}\sup_{\Lambda\geq\lambda}\sup_{\alpha \in I_M}\left
    \vert\frac{\sum_{k=\lambda}^{\Lambda}a_k\left[\exp{2i\pi\alpha j X_k(\omega)}-{\EE}\exp{2i\pi\alpha j X_k}\right]}
    {\sqrt{A_{\lambda,\Lambda,M}^2\log{(\vert j \vert +3)}}}\right \vert < \infty$$
\end{theo}

\begin{rem}
When$(X_k)_{k\geq 0}$ takes integer values, the proof of
theorem \ref{principal} is easier. Namely, using the fact that $\alpha\mapsto
j\alpha$ (mod 1) is onto for $j\not=0$, we get
  $$\sup_{j\in{\Z}^*}\sup_{\alpha \in {\T}}\left
   \vert \sum_{k=\lambda}^{\Lambda} a_k\left[\exp {2i\pi \alpha
   j X_k(\omega) } - {\EE}\exp {2i\pi \alpha  j X_k } \right ] \, \right \vert$$
  $$=\sup_{\alpha \in {\T}}\left \vert \sum_{k=\lambda}^{\Lambda}  a_k \left [\exp {2i\pi \alpha
   X_k(\omega) } - {\EE}\exp {2i\pi \alpha   X_k }\right ] \,\right \vert$$

the result of theorem \ref{principal} becomes then :
  $${\EE}\sup_{\lambda\geq 1}\sup_{\Lambda\geq\lambda}\sup_{\alpha \in I_M}\left
    \vert\frac {\sum_{k=\lambda}^{\Lambda}a_k\left[\exp {2i\pi\alpha X_k(\omega)}-{\EE}\exp{2i\pi\alpha X_k}
    \right]}{\sqrt{A_{\lambda,\Lambda,M}^2}}\right \vert<\infty $$
When $(X_k)_{k\geq 1}$ takes real values, the proof is more tedious. It relies on a fine inequality
about decoupling gaussian random functions (see section 3.). We can see here why, for integer-valued $X_k$, we
can work with the functional space $A({\T})$, whereas for real-valued $X_k$, we need to introduce the space 
$B({\T})$.
\end{rem}

\section{Proof of theorem \ref{moment} and corollary \ref{separation}}

First, we split $F$ into two parts as follows :
  $$\sum_k a_k f(\alpha X_k(\omega))=\sum_k a_k[f(\alpha X_k(\omega)-{\EE}(f(\alpha X_k))
    +\sum_k a_k{\EE}(f(\alpha X_k))$$
\par
{\bf -Step 1 : }(first part of the sum)\\
Let $(N_k)_{k\geq 1}$ be a strictly increasing sequence of integers and define
  $$P_k(\alpha)=\sum_{l= N_k +1}^{N_{k+1}}a_l \left [f{(\alpha X_l(\omega))}
    -{\EE} f{(\alpha X_l )}\right]\leqno { \forall k\geq 1,}$$
where $f\in B({\T})$.
We want to study the following series, for all $M\geq 1$ :
  $$\sum_k\sup_{\alpha\in [-M, M]}\vert P_k(\alpha)\vert$$
We have :
  $$\vert P_k(\alpha) \vert \leq\sum_{j\in {\Z}}\vert\hat{f}(j)\vert\left\vert
    \sum_{l= N_k +1}^{N_{k+1}}a_l[\exp{(2\pi j\alpha X_l(\omega))}-{\EE} \exp{(2\pi j\alpha X_l)}]
    \right \vert$$
Hence, using theorem \ref{principal}, there exists a positive integrable random variable
$\xi$ such that
\begin{equation}\label{majoration}
  \sup_{\alpha \in[-M,M]}\vert P_k(\alpha)\vert\leq\xi\vert\vert\vert f\vert\vert\vert
  \sqrt{\log{(M\Phi_{\beta}(N_{k+1}))}\sum_{j=N_k +1}^{N_{k+1}}\vert a_j\vert ^2}
\end{equation}
where
  $$\xi=\sup_{j\in {\Z}}\sup_{k\geq 1}\sup_{\alpha \in I_M}\left
    \vert\frac{1}{\sqrt{A_{k,M}^2\log{(\vert j\vert +3)}}}\sum_{l=N_k +1}^{N_{k+1}}a_l\left[e^{2i\pi\alpha
    j X_l(\omega)}-{\EE}e^{2i\pi\alpha jX_l}\right]\right\vert$$
with
  $$A_{k,M}^2={\log{(M\Phi_{\beta}(N_{k+1}))}\sum_{l= N_k +1}^{N_{k+1}}\vert a_l \vert ^2}$$
First, in the polynomial case, that is to say when there exists $d>0$ with
  $$\Phi_{\beta}(N)={\mathcal O}(N^d)$$
then we choose $N_k=2^{2^k}$ and we need to prove that
  $$\sum_k2^{k/2}\left(\sum_{l=2^{2^k}+1}^{2^{2^{k+1}}}\vert a_l \vert ^2\right)^{1/2}<+\infty$$
now we use the following equivalent:
  $$\sum_{l=2^{2^k}+1}^{2^{2^{k+1}}}\frac{1}{l(\log(l))^{1/2}}\approx 2^{k/2}$$
which may be computed by comparing series and integral, hence :
\begin{eqnarray*}
  2^{k/2}\left(\sum_{l=2^{2^k}+1}^{2^{2^{k+1}}}\vert a_l \vert ^2\right)^{1/2} & \le & 
  C\sum_{l=2^{2^k}+1}^{2^{2^{k+1}}}\frac{1}{l(\log(l))^{1/2}}\left(\sum_{l=2^{2^k}+1}^{\infty}\vert a_l \vert ^2\right)^{1/2} \\
  & \le &C\sum_{l=2^{2^k}+1}^{2^{2^{k+1}}}\frac{\left(\sum_{j=l}^{\infty}\vert a_j \vert ^2\right)^{1/2}}{l(\log(l))^{1/2}}
\end{eqnarray*}
and, using condition \ref{condition1} :
\begin{eqnarray*} 
  \sum_k2^{k/2}\left(\sum_{l=2^{2^k}+1}^{2^{2^{k+1}}}\vert a_l \vert ^2\right)^{1/2} & \le & 
  \sum_k C\sum_{l=2^{2^k}+1}^{2^{2^{k+1}}}\frac{\left(\sum_{j=l}^{\infty}\vert a_j \vert ^2\right)^{1/2}}{l(\log(l))^{1/2}} \\
  & \le & \sum_{n\geq 2}\frac{\sqrt{\sum_{k\geq n}\vert a_k\vert^2}}{n\sqrt{\log n}}<+\infty
\end{eqnarray*}
this implies :
  $$\sum_{k\geq 1}\sup_{\alpha \in [-M, M]}\vert P_k(\alpha) \vert < \infty$$
almost everywhere on the measurable set $\Omega_o =\{\omega\in\Omega, \xi (\omega)<\infty\}$.
By construction, this set does not depend on the choice of $f$. \\
Secondly, in the subexponential case, that is when there exists $\gamma\in]0,1[$ with
  $$\Phi_{\beta}(N)={\mathcal O}(2^{N^{\gamma}})$$
we choose $N_k=2^k$ and we need to prove that
  $$\sum_k2^{\gamma k/2}\left(\sum_{l=2^k+1}^{2^{k+1}}\vert a_l \vert ^2\right)^{1/2}<+\infty$$
Using the following equivalent:
  $$\sum_{l=2^k+1}^{2^{k+1}}\frac{1}{l^{1-\frac{\gamma}{2}}}\approx 2^{\gamma k/2}$$
and doing the same kind of computation as before, using condition \ref{condition2}:
  $$ \sum_{n\geq 2}\frac{\sqrt{\sum_{k\geq n}\vert a_k\vert^2}}{n^{1-\frac{\gamma}{2}}}<+\infty$$
implies
$$\sum_{k\geq 1}\sup_{\alpha\in[-M, M]}\vert P_k(\alpha)\vert<\infty$$
We also get from (\ref{majoration}), for all $\alpha\in{\R}$:
  $$\vert P_k(\alpha)\vert\le\xi\vert\vert\vert f\vert\vert\vert\sqrt{\log(\vert\alpha\vert+1)}
  \sqrt{\log(\Phi_{\beta}(N_{k+1}))\sum_{j=N_k +1}^{N_{k+1}}\vert a_j\vert ^2}$$
summing on $k$, we get the inequality:
  $$\vert F(\alpha,\omega)-{\EE}(F(\alpha,.))\vert\le C\xi(\omega)\vert\vert\vert f\vert\vert\vert\sqrt{\log(\vert \alpha\vert+2)}$$
where $C$ only depends on $(a_k)$ and ${\EE}(\xi)<\infty$\\
This ends the proof of theorem \ref{centered} which
is also the first step of the proof of theorem \ref{moment}\\

\par

{\bf -Step 2 : }(second part of the sum)\\
Let $K$ be a compact which does not contain zero and $\alpha\in K$. Let $n<m$ be two integers,
\begin{eqnarray}
  \vert\sum_{k=n}^m a_k{\EE}(f(\alpha X_k))\vert &=& \vert\sum_{j\in{\Z}^*}\sum_{k=n}^m a_k\hat f(j)
  {\EE}(\exp(2i\pi j\alpha X_k))\vert\\
  &=& \vert\sum_{j\in{\Z}^*}\sum_{k=n}^m a_k\hat f(j)\varphi_{X_k}(j\alpha)\vert\\
  &=& \vert\sum_{j\in{\Z}^*}\hat f(j)(\sum_{k=n}^m a_k\varphi_{X_k}(j\alpha))\vert\\
  &\le& (\sum_{j\in{\Z}^*}\vert\hat f(j)\vert)\sup_{j\in{\Z}^*}\vert\sum_{k=n}^m a_k\varphi_{X_k}(j\alpha)\vert
  \label{abel}
\end{eqnarray}
At this point, to prove theorem \ref{moment}, we can conclude directly by using hypothesis $({\mathcal H})$ to get
$$\sup_{n<m}\sup_{\alpha\in K}\vert\sum_{k=n}^ma_k{\EE}(f(\alpha X_k))\vert<\varepsilon$$
as long as $m$ and $n$ are large enough.\\
To prove corollary \ref{separation}, we use Abel's summation. Let $\phi_p=\sum_{k=0}^{p-1}\varphi_{X_k}$, we have :
\begin{eqnarray}
  \sum_{k=n}^ma_k\varphi_{X_k}(j\alpha)&=&\sum_{k=n}^ma_k(\phi_{k+1}(j\alpha)-\phi_k(j\alpha))\\
  &=&\sum_{k=n+1}^{m+1}a_{k-1}\phi_k(j\alpha)-\sum_{k=n}^ma_k\phi_k(j\alpha)\\
  &=&-a_n\phi_n(j\alpha)+a_m\phi_{m+1}(j\alpha) \\
  &+&\sum_{k=n+1}^{m}(a_{k-1}-a_k)\phi_k(j\alpha)
\end{eqnarray}
and :
  $$(\ref{abel})\le\vert\vert\vert f\vert\vert\vert\sup_{N\geq 1}\sup_{\alpha\in K}\sup_{j\in{\Z}^*}
  \vert\phi_N(j\alpha)\vert\left\lbrack\vert a_m\vert+\vert a_n\vert+\sum_{k=n+1}^m\vert a_k-a_{k-1}\vert\right\rbrack$$
we now conclude using hypothesis ${\mathcal H}'$ and hypothesis (1) on the sequence $(a_n)$ in corollary \ref{separation} :
  $$\sup_{n<m}\sup_{\alpha\in K}\vert\sum_{k=n}^ma_k{\EE}(f(\alpha X_k))\vert$$
  $$\le\vert\vert\vert f
  \vert\vert\vert\left\lbrack\vert a_m\vert+\vert a_n\vert+\sum_{k=n+1}^{m}\vert a_k-a_{k-1}\vert\right\rbrack
  \sup_{N\geq 1}\sup_{\alpha\in K}\sup_{j\in{\Z}^*}\vert\phi_N(j\alpha)\vert<\varepsilon$$
as long as $m$ and $n$ are large enough.\\
To prove corollary \ref{separationcauchy}, we use Cauchy Schwarz inequality in the following way:
  $$\left\vert\sum_{k=n}^m a_k\varphi_{X_k}(j\alpha)\right\vert\le\sqrt{\sum_{k=n}^m\vert a_k\vert^2}
  \sqrt{\sum_{k=n}^m\vert\varphi_{X_k}(j\alpha)\vert^2}$$
and we conclude using condition $H''$.

\begin{rem}\label{remarque}\ \\
\begin{enumerate}
\vspace{-0.5cm}
\item the most general hypothesis we can put on the sequence $(a_k)_{k\geq1}$ is the following :
there exists a strictly increasing sequence 
$(N_k)_{k\geq 1}$ such that
  $$\sum_{k=1}^{\infty}\sqrt{\log\Phi_{\beta}(N_{k+1})\sum_{l=N_k +1}^{N_{k+1}}\vert a_l\vert ^2}<\infty$$
\item if the process $(X_k)_{k\geq 1}$ takes integer values then $f\in A({\T})$ can be assumed
without any other hypothesis.
\item If $\sum \vert a_k\vert^2$ diverges, then we can construct a stochastic
process  $(X_k)_{k\geq 1}$ verifying the hypothesis of theorem \ref{centered} and find $f\in B({\T})$ 
such that the convergence of the series is not uniform on any compact.
In that sense, the conditions imposed to the sequence
$(a_k)_{k\geq 1}$ are optimal. Remark that in this case, condition \ref{condition1} is not fulfilled. \par
Namely consider a sequence of
independent random variables with disjoint supports.
For all $k\geq 1$ the support of $X_k$ is the set of integers belonging to
$[k^2,(k+1)^2-1]$ and hence, the hypothesis on the moment is verified. We will come back to the law
of $X_k$ later. Now choose $f$ in the following way : for all $\alpha\in {\T}$,
$f(\alpha)=\exp{(2i\pi\alpha)}.$ Thus $f\in A({\T})$  et $\vert\vert\vert f\vert\vert\vert <\infty$.
As a consequence, if the convergence of the series defining $F$ was uniform in
$\alpha$ on ${\T}$, then we would have :
\begin{eqnarray*}
  &&\sqrt{\int_0^1\left\vert\sum_{k=1}^{\infty}a_k\left[\exp{2i\pi\alpha
    X_k(\omega)}-{\EE}\exp{2i\pi\alpha X_k}\right]\right\vert^2 d\alpha}\\
  &\leq& \sup_{\alpha\in[0,1]}\left\vert\sum_{k=1}^{\infty}a_k\left[\exp{2i\pi\alpha
    X_k(\omega)}-{\EE}\exp{2i\pi\alpha X_k}\right]\right\vert<\infty 
\end{eqnarray*}
By construction, ${\PP}-$ almost surely :
\begin{eqnarray*}
  &&\int_0^1\left\vert\sum_{k=1}^{\infty}a_k\left[\exp{2i\pi\alpha
  X_k(\omega)}-{\EE}\exp{2i\pi\alpha X_k}\right]\right\vert^2d\alpha\\
  &=&\sum_{k=1}^{\infty}\vert a_k\vert^2\int_0^1\vert\exp{(2i\pi\alpha X_k(\omega))}-{\EE}\exp
  {(2i\pi\alpha X_k)}\vert^2 d\alpha \\
  &\geq&\sum_{k=1}^{\infty}\vert a_k\vert^2\int_0^1\left\vert 1-\vert{\EE}\exp{(2i\pi\alpha X_k)}\vert
  \right\vert^2 d\alpha
\end{eqnarray*}
Assume now that  the law of $X_k$ is uniform on the $2k +1$ integers of
$[k^2, (k+1)^2-1].$ for all $k\geq 1$.
  $$\vert{\EE}\exp{(2i\pi\alpha X_k)}\vert=\frac{1}{2k+1}\left\vert\frac{\sin{\pi\alpha(2k+1)}}{\sin{\pi\alpha}}\right\vert$$
Using Lebesgue convergence theorem, we get
  $$\lim_{k\rightarrow +\infty}\int_0^1\left\vert 1-\vert{\EE}\exp{(2i\pi\alpha X_k)}\vert\right\vert^2 d\alpha=1$$
and we also get the divergence of the series with positive terms
  $$\sum_{k=1}^{\infty}\vert a_k\vert^2\int_0^1\left\vert 1-\vert{\EE}\exp{(2i\pi\alpha X_k)}\vert\right\vert^2 d\alpha=\infty$$
A contradiction with uniform convergence of the centered part.
\end{enumerate}
\end{rem}

\section{Proof of theorem \ref{principal}}

Let us begin by restating some inequalities obtained by Fernique (\cite{Fernique})
which will be useful in the proof of theorem \ref{principal}.
For more information on gaussian techniques in this framework, see \cite{SchneiderHabi} and \cite{DurSchnei}.

\begin{ineq}\label{ineg1}
Let $( G_k)_{ k\geq1}$ be a sequence of Banach space valued gaussian random variables
$(B,\parallel\cdot\parallel)$ defined on a probabilised space $(\Omega,
\mathcal A,\mathbb P)$. Then :
  $${\EE}\sup_{k\geq 1}\parallel G_k\parallel\leq K_1\left\{\sup_{k\geq 1}{\EE}\parallel G_k\parallel+{\EE}\sup_{k\geq 1}
    \mid\lambda _k\sigma_k\mid\right\}$$
where $( \lambda _k)_{ k\geq1 }$ is an isonormal sequence, $K_1$ a universal constant
and for all $k\geq 1$,
  $$\sigma _k =\sup_{f\in B^{'},\parallel f\parallel\leq 1}\parallel<G_k,f>_{B}\parallel_{2,\mathbb P}$$
\end{ineq}

\begin{ineq}\label{ineg2}
Let $g$ be a real valued stationary gaussian random variable, separable and continuous in quadratic mean.
Let $m$ be its associated spectral measure on ${\mathbb {R}}^{+}$ defined by
  $${\EE}[\mid g(s)-g(t)\mid ^2]=2\int_0^\infty[1-\cos{2\pi u(s-t)}]m(du)$$
We have 
  $${\EE}\sup_{\alpha\in[0,1]}g(\alpha)\leq K\left\{\sqrt{\int_0^\infty\min(u^2,1)m(du)}+
    \int_0^\infty\sqrt{m\left(\left]e^{x^2},\infty\right[\right)}dx\right\}$$
where $K$ is a universal constant.
\end{ineq}

\begin{ineq}\label{ineg3}({\bf decoupling})
Let $X=\{X(t):t\in {\T}\}$ be a gaussian random function defined on a finite or countable set
$T$. Let $\{T_k,k\in[1,n]\}$ be a covering of $T$.
Let $S = T_1 \times \cdots \times T_n$. The following inequality holds :
  $${\EE}\left\{\sup_{k\in[1,n]}\left[\sup_{t\in T_k}X(t)-{\EE}\sup_{t\in T_k}X(t)
    \right]\right\}\leq\frac{\pi}{2}\cdot\sup_{s\in S}\left\{{\EE}\left[\sup_{k\in[1,n]}X(s_k)\right]\right\}$$
\end{ineq}

The following estimation generalizes inequality \ref{ineg2} and will be useful in the almost periodic case because
it gives estimations on arbitrarily large intervals.

\begin{ineq}\label{ineg4}
Let $g$ a real valued stationnary gaussian random function, separable and continuous in quadratic mean.
Let $m$ its associated spectral measure on ${\mathbb {R}}^{+}$ defined as in inequality \ref{ineg2}.
There exists a universal constant $K$ such that
  $${\EE}\sup_{\alpha\in[-M,M]}g(\alpha)$$
  $$\leq K\left\{\sqrt{\int_0^\infty\min(2Mu^2,1)\ m(du)}
    +\int_0^\infty\sqrt{m\left(\left]\frac{e^{x^2}}{2M},\infty\right[\right)}dx \right \}$$
\end{ineq}

Let us now come to the proof

{\bf -Step 1:} 
In this part, we replace our problem by a question of regularity of trajectories of random gaussian functions.
\par
Let us consider an independent copy of  $X=(X_k)_{k\geq 1}$ denoted by
$X'=(X_k')_{k\geq 1}$ defined on another probabilized space $(\Omega^{'}, {\mathcal A}^{'} , {\mathbb P}^{'})$. We call
${\EE}_{*}$ the integration symbol whose index refers to the space of integration.
\par
Using classical convexity properties, to prove (1), it is enough to show
\begin{eqnarray}\label{convexity}
  {\EE}_{X,X^{'}}\sup_{j\in{\Z}}\sup_{\lambda\geq 1}\sup_{\Lambda\geq\lambda}\sup_{\alpha \in I_M}\left\vert
    \frac{\sum_{k=\lambda}^{\Lambda}a_k\left[e^{2i\pi\alpha jX_k}-e^{2i\pi\alpha jX^{'}_k}\right]}
    {\sqrt{A_{\Lambda,\lambda,M}^2\log{(\vert j\vert +3)}}}\right\vert<\infty
\end{eqnarray}
\par

Let us now symmetrize the problem: consider the following separable family of random functions, with continuous
trajectories
  $$(f_k)_{k\geq 1}=\{f_{k}(\alpha, j)=a_k(\exp 2i\pi\alpha jX_k- \exp 2i\pi\alpha jX'_k),\alpha\in I_M,j\in{\Z}\}_{ k\geq1}$$
By construction  $f$ is a symmetric family of random functions, that is to say their law is sign-invariant. More precisely, call
$\{\varepsilon_{k},k\geq 1\}$ a sequence of independent Rademacher random variables
(taking the values $+1$ and $-1$ with probability $1/2$), defined on a third space $(\Omega^{''}, {\mathcal A}^{''} , {\mathbb P}^{''})$,
independent of $X$ and $X'$. $\{f_{k}, k\geq 1\}$ and $\{\varepsilon_{k} f_{k}, k\geq 1\}$ have the same law.
Thus for all integers $(\Lambda,\lambda)$ such that $\Lambda\geq\lambda$, $\sum_{k=\lambda}^{\Lambda}f_{k}$ and
$\sum_{k=\lambda}^{\Lambda}\varepsilon_{k}f_{k}$ also have the same law.\par
That is why (\ref{convexity}) can be written on a larger space of integration in the following way:
  $${\EE}_{X,X^{'},\varepsilon}\sup_{j\in{\Z}} \sup_{\lambda\geq 1}\sup_{\Lambda\geq\lambda}\sup_{\alpha \in I_M}\left
    \vert\frac{\sum_{k=\lambda}^{\Lambda }\varepsilon_{k}f_{k}(\alpha,j)}{\sqrt{A_{\Lambda,\lambda,M}^2\log{(\vert j\vert +3)}}}
    \right\vert$$
We deduce a sufficient condition for (\ref{convexity}) to be realised
\begin{eqnarray}\label{epsilon}
  {\EE}_{X,\varepsilon}\sup_{j\in{\Z}}\sup_{\lambda\geq 1}\sup_{\Lambda\geq\lambda}\sup_{\alpha \in I_M}\left
   \vert\frac{\sum_{k=\lambda}^{\Lambda }\varepsilon_k a_k\exp{2i\pi\alpha jX_k}}{\sqrt{A_{\Lambda,\lambda,M}^2\log{(\vert j\vert +3)}}}
   \right\vert
\end{eqnarray}
We then use a precious tool in the theory of gaussian random functions: the contraction principle. This tool is built on a quite
simple idea: replace the choice of signs by a sequence of gaussian random variables with mean zero and variance one. This idea can be 
explained by the following property: given $g$ a gaussian random variable  with mean zero and variance 1 and 
$\varepsilon$  a Rademacher random variable, if $g$ and $\varepsilon$ are independent, then $g$ and
$\varepsilon |g|$ have the same law.

\par
As a consequence, in order to prove (\ref{epsilon}), we show
  $${\EE}_{X,g,g^{'}}\sup_{j\in{\Z}}\sup_{\lambda\geq 1}\sup_{\Lambda\geq\lambda}\sup_{\alpha\in I_M}\left
   \vert\frac{\sum_{k=\lambda}^{\Lambda}a_k(g_{k}\cos 2\pi\alpha jX_k+g'_{k}\sin 2\pi\alpha jX_k)}{\sqrt
   {A_{\lambda,\Lambda,M}^2\log{(\vert j\vert +3)}}}\right\vert<+\infty$$
where $\{g_{k},k\geq 1\}$ et $\{g'_{k},k\geq 1\}$ are two sequences of 
independent identically distributed random variables with law
${\mathcal N}(0,1)$, independent of $X$ and $\varepsilon$, defined on two other probabilised spaces. \par
Conditionally to $X$, the problem is reduced to studying the regularity of the trajectories of stationary gaussian
random variables. This concludes the first step of the proof.
\medskip\medskip
\par

{\bf -Step 2 :} In this part, we use the gaussian tools introduced in the beginning.

Conditionally to $X$, call $G(\lambda,\Lambda,j,\alpha)$ the following quantity
  $$\frac {1}{\sqrt{A_{\lambda,\Lambda,M}^2\log{(\vert j\vert +3)}}}
    \sum_{k=\lambda}^{\Lambda}a_k\left[g_{k}\cos{(2\pi\alpha jX_k)}+g^{'}_{k}\sin{(2\pi\alpha jX_k)}\right]$$
If $j$, $\lambda$ and $\Lambda$ are fixed, $G(\alpha):=G(\lambda,\Lambda,j,\alpha)$ is a random function with almost surely
continuous trajectories (up to a modification of trajectories). That is why it is enough to show that $G$ is bounded on
$I_M\cap{\mathbb Q}$. Moreover, we will assume that $\vert j\vert\leq J$ where $J$ is a large fixed integer.

\par
Let us begin by finding an upper bound for
  $${\EE}_{g,g^{'}}\sup_{\vert j\vert\leq J}\sup_{\lambda\geq 1}\sup_{\Lambda\geq \lambda}
    \sup_{\alpha\in I_M\cap{\mathbb Q}}\vert G(\lambda,\Lambda,j,\alpha)\vert$$
First remark that if $Y_t$ $(t\in T)$ is a gaussian random function defined on $T$, then for all
$t_0\in T$ we have (see \cite{Fernique} (480))
  $${\EE}\sup_{t\in T}\vert Y_t\vert\leq{\EE}\vert Y_{t_0}\vert+{\EE}\sup_{t\in T} Y_t$$
In this way, we get rid of the absolute value. Apply this remark to
$G(\lambda,\Lambda,j,\alpha)$ with $\alpha=0$, $\lambda=\Lambda=1$ and $j=0$
and let us find an upper bound for
\begin{eqnarray}
  {\EE}_{g,g^{'}}\sup_{\vert j \vert \leq J}\sup_{\lambda\geq 1}\sup_{\Lambda\geq\lambda}
    \sup_{\alpha\in I_M\cap{\mathbb Q}}G(\lambda,\Lambda,j,\alpha)
\end{eqnarray}
In order to apply the decoupling inequality \ref{ineg3},  define
  $$T=\{-J,\cdots,J\}\times{H}\times(I_M\cap{\mathbb Q})$$
where $J$ is a large enough integer and $H$ is the upper triangle of 
dimension 2 in ${\N}\times{\N}$ (see figure 1 below) ($\lambda\in{\N}$ and $\Lambda \geq \lambda$.)
A point in $t\in T$ will be written $t=(j,\lambda,\Lambda,\alpha)$.
\par
This set $T$ is at most countable and we will find an upper bound for
${{\EE}_{g, g^{'}}}\sup_{t\in T}G(t)$ independently of $J$ and then
conclude by taking the supremum on $j\in J$.
Define
  $$T_j=\{ j\}\times{H}\times(I_M\cap{\mathbb Q})$$
It is obvious that $\{T_j\}_{j=-J,\cdots,J}$ is a covering of $T$.
Define 
  $$S=T_{-J}\times\cdots\times T_J$$
Using inequality \ref{ineg3}, we have
  $${\EE}_{g,g^{'}}\sup_{t\in T}G(t)\leq\frac{\pi}{2}\sup_{s\in S}{{\EE}_{g,g^{'}}}\sup_{-J\le j\le J}G(s_j)
    +\sup_{-J\le j\le J}{{\EE}_{g,g^{'}}}\sup_{t\in T_j}G(t)$$
where $s_j$ is a point in $T_j.$
\medskip \medskip

{\bf -Step 3 :} We study now
  $$\sup_{s\in S}{{\EE}}\sup_{-J\le j\le J}G(s_j).$$
We can rewrite this in the following way :
\begin{eqnarray}\label{j}
  \sup{{\EE}_{g,g^{'}}}\sup_{-J\le j\le J}\frac{\sum_{k=\lambda_{j}}^{\Lambda_{j}}a_k\left[g_{k}\cos{(2\pi j\alpha_j X_k )}
   +g^{'}_{k}\sin{(2\pi j\alpha_j X_k)}\right]}{\sqrt{A_{\lambda_j,\Lambda_j,M}^2\log{(\vert j\vert +3)}}}
\end{eqnarray}
where the first supremum is taken on
$$\{(\alpha_{-J},\cdots,\alpha_J)\in(I_M\cap{\mathbb Q})^{2J+1},(\, 
(\lambda_{-J},\Lambda_{-J})\, , \cdots ,
(\lambda_{J},\Lambda_{J}) \, )\in {H}^{2J+1}\,\}$$
\par
Fix $ \{(\alpha_{-J},\cdots,\alpha_J)\in(I_M\cap{\mathbb Q})^{2J+1}$ and $((\lambda_{-J},\Lambda_{-J})\, , \cdots ,
(\lambda_{J},\Lambda_{J}) \, )\in {H}^{2J+1}\,\}$. Define the gaussian process
  $$G_{j}:=\frac{\sum_{k=\lambda_{j}}^{\Lambda_{j}}a_k\left[g_{k}\cos{(2\pi j\alpha_j X_k )}+g^{'}_{k}
    \sin{(2\pi j\alpha_j X_k)}\right]}{\sqrt{A_{\lambda_j,\Lambda_j,M}^2\log{(\vert j\vert +3)}}}$$
In order to get an upper bound for \ref{j}, we remark that
$$\sup{\EE}_{g,g^{'}}\sup_{-J\le j\le J}G_{j}$$
\begin{eqnarray}\label{supesp}
  \leq\sup_{((\lambda_{-J},\lambda_{-J}),\cdots ,
    (\lambda_{J},\lambda_{J}))\in{H}^{2J+1}}{\EE}_{g,g^{'}}\sup_{-J\le j\le J}\sup_{\alpha \in I_M}
    \vert G(j,\lambda_j,\Lambda_j,\alpha)\vert
\end{eqnarray}
We will get an upper bound for the right hand side of \ref{supesp}
independently of $J$, which will give us an upper bound for \ref{j} by taking the supremum on $j\in J$.
\par

Applying inequality  \ref{ineg1} to the finite sequence of random gaussian functions
  $$(G(j,\lambda_j,\Lambda_j,\alpha))_{-J\leq j\leq J}$$
we prove that $\mathbb  E_{g,g^{'}}\sup_{-J\leq j\leq J}\sup_{\alpha\in I_M}\vert G(j,\lambda_j,\Lambda_j,\alpha)\vert$
is less than
\begin{eqnarray}\label{isonormal}
  C\left\{\sup_{-J\leq j\leq J}{\EE}_{g,g^{'}}\sup_{\alpha\in I_M}
    \vert G(j,\lambda_j,\Lambda_j,\alpha)\vert+{\EE}_{\xi}\sup_{-J\leq j\leq J}\vert\xi_{j}q_{j}\vert\right\}
\end{eqnarray}
where  $C$ is a universal constant, $(\xi_j)_{-J\leq j\leq J}$ is an isonormal sequence and
  $$q_{j}\leq\sup_{\alpha\in I_M}\left\vert\left\vert G(j,\lambda_j\Lambda_j,\alpha)\right\vert\right\vert_{2,g,g^{'}}$$
This gives us the following upper bound
  $$q_{j}\leq \frac {1}{\sqrt{\log(M\Phi_{\beta}(\Lambda_j))\log{(\vert j\vert +3)}}}$$
As for all $j$ we have $\Lambda_j\geq 1$ and $M\geq 1$ we easily get
\begin{eqnarray}
  \mathbb  E_{\xi}\sup_{-J\leq j\leq J}\vert\xi_{j} q_{j}\vert&\leq&{\EE}_{\xi}\sup_{-J\leq j\leq J}\left\vert
  \frac {1}{\sqrt{\log{(\vert j\vert +3)}}}\xi_{j}\right\vert\\
  &\leq& C{\EE}_{\xi}\sup_{j\in{\Z}}\left\vert\frac{\xi_{j}}
    {\sqrt{\log(\vert j\vert +3)}}\right\vert
\end{eqnarray}
The exponential integrability of gaussian vectors gives
  $${\EE}_{\xi}\sup_{j\in{\Z}}\left\vert\frac{\xi_{j}}
    {\sqrt{\log{{{\vert j\vert +3}}}}}\right\vert=C<\infty$$
and hence 
  $$\sup_{(\lambda_{-J},\Lambda_{-J}),\cdots,(\lambda_{J},\Lambda_{J}))\in {H}^{2J+1}}
    {\EE}_{\xi}\sup_{-J\leq j\leq J}\vert\xi_jq_j\vert=C<\infty$$
independently of $J$ and of the sequence $(\lambda_{-J},\Lambda_{-J}),\cdots,(\lambda_{J},\Lambda_{J}))\in {H}^{2J+1})$.
\par
We then get independently of $J$ on the whole integration space
  $${\EE}_{X}\sup_{(\lambda_{-J},\Lambda_{-J}),\cdots,(\lambda_{J},\Lambda_{J}))\in {H}^{2J+1}}{\EE}_{\xi}\sup_{-J\leq j\leq J}
  \vert\xi_{j}q_{j}\vert\leq{C}<\infty$$
where $C$ is a constant.

\medskip \medskip
\par

Let us now try to find an upper bound for
\begin{eqnarray}\label{supg}
  \sup_{((\lambda_{-J},\Lambda_{-J}),\cdots,(\lambda_{J},\Lambda_{J}))\in{H}^{2J+1}}
    \sup_{-J\leq j\leq J}{\EE}_{g,g^{'}}\sup_{\alpha \in I_M}\vert G(j,\lambda_j,\Lambda_j,\alpha)\vert
\end{eqnarray}
independently of $J$.

\par
We will use inequality  \ref{ineg4}.
Let us choose a finite sequence 
$$((\lambda_{-J},\Lambda_{-J}),\cdots,(\lambda_{J},\Lambda_{J})$$
and an integer
$\vert j\vert\leq J$. The gaussian random function $G_{j}(\alpha):=G(j,\lambda_j,\Lambda_j,\alpha)$ is stationary.
Its associated spectral measure on ${\R}^+$ is defined by
  $$m_{j}=\frac {1}{\log{(\vert j\vert +3)}\log{(M\Phi_{\beta}(\Lambda_j))}\left(\sum_{k=\lambda_j}^{\Lambda_j}\vert a_k\vert ^2\right)}
    \sum_{k=\lambda_j}^{\Lambda_j}\vert a_k\vert^2\delta _{j\vert X_k\vert}$$
where $\delta_u$ is the Dirac measure in the point $u$.
\par

We get
  $${\EE}_{g,g^{'}}\sup_{\alpha\in I_M}\vert G_{j}(\alpha)\vert$$
    $$\leq C\left\{\sqrt{\int_0^\infty\min(\frac{u^2}{2M},1)m_{j}(du)}+
    \int_0^\infty\sqrt{m_{j}\left(]\frac{e^{x^2}}{2M},\infty[\right)}dx\right\}$$
It is obvious that the first term is less than 
  $$\sqrt{m_{j}({\R}^+)}\leq C{\frac{1}{\sqrt{\log{(2M)}\log{(\vert j\vert +3)}}}}$$
where $C$ is a universal constant because for all $j$ we have $\Phi_{\beta}(\Lambda_j)\geq 2$.
\par

For the second term, it can be rewritten in the following way
\begin{eqnarray*}
  &&\int_0^\infty \sqrt{m_{j}\left(]\frac{\exp{x^2}}{2M},\infty[\right)}dx\\
  &=&\left[\frac{1}{\log{(\vert j\vert +3)}\log{(M\Phi_{\beta}(\Lambda_j))}\left(\sum_{k=\lambda_j}^{\Lambda_j}\vert a_k\vert ^2\right)}
    \right]^{\frac{1}{2}}\\
  &\times&\int_0^\infty\sqrt{\sum_{k=\lambda_j}^{\Lambda_j}\vert a_k\vert^2{\bf 1}_{\{2M\mid jX_k\mid>e^{x^2}\}}}dx
\end{eqnarray*}
Using
  $${\bf 1}_{\{2M\mid jX_k\mid>e^{x^2}\}}\leq{\bf 1}_{\{2M\vert j\vert\sup_{l\leq k}\mid X_l\mid>e^{x^2}\}}\leqno {\forall k \geq 1,}$$
we cut ${\mathbb  R}^+$ in the integral according to the increasing subdivision
  $$\{0\}\bigcup\{\sqrt{\log^{+}{(2M\vert j\vert\sup_{l\leq k}\mid X_l\mid)}},\lambda_j\leq k\leq\Lambda_j\}$$

We get thus an upper bound for the previous integral
\begin{eqnarray*}
\sum_{k=\lambda_j}^{\Lambda_j}&&\left[\sqrt{\log^{+}{(2M\vert j\vert\sup_{l\leq k}\mid X_l\mid)}}-\sqrt{\log^{+}
    {(2M\vert j\vert\sup_{l\leq k-1}\mid X_l\mid)}}\right]
    \left(\sum_{l=\lambda_j}^{\Lambda_j}\vert a_l\vert^2\right)^{\frac{1}{2}}\\
  &\leq&\left(\sum_{l=\lambda_j}^{\Lambda_j}\vert a_l\vert^2\right)^{\frac{1}{2}}\sum_{k=\lambda_j}^{\Lambda_j}\left[
    \sqrt{\log^{+}{(2M\vert j\vert\sup_{l\leq k}\mid X_l\mid)}}\right.\\
    &-&\left.\sqrt {\log^{+}{(2M\vert j\vert\sup_{l\leq k-1}\mid X_l\mid)}}\right]\\
  &\leq& 2\left(\sum_{k=\lambda_j}^{\Lambda_j}\vert a_k\vert^2\right)^{\frac{1}{2}}
    \sqrt{\log^{+}{(2M\vert j\vert\sup_{l\leq\Lambda_j}\mid X_l\mid)}}
\end{eqnarray*}
Consequently, $\forall\vert j\vert\leq J$,
\begin{eqnarray*}
  \int_0^\infty\sqrt{m_{j}\left(]\frac{\exp{x^2}}{M},\infty[\right)}dx&\leq& 2\sqrt{\frac{\log^{+}
    {(2M\vert j\vert\sup_{l\leq\Lambda_j}\mid X_l\mid)}}{\log{(\vert j\vert +3)}\log(M\Phi_{\beta}(\Lambda_j))}}\\
  &\leq& 4\sqrt{\frac{\log^{+}{(2M\sup_{1\leq l\leq\Lambda_j}\mid X_l\mid)}}{\log(M\Phi_{\beta}(\Lambda_j))}}\\
  &\leq& 4\sup_{1\leq l\leq\Lambda_j}\sqrt{\frac{\log^{+}{(2M\mid X_l\mid)}}{\log(M\Phi_{\beta}(l))}}\\
  &\leq& 4\sup_{N\geq 1}\sqrt{\frac{\log^{+}{(2M\mid X_N\mid)}}{\log{(M\Phi_{\beta}(N))}}}\\
  &\leq& 4\left(1+\sup_{N\geq 1}\sqrt{\frac{\log^{+}
    {(\mid X_N\mid)}}{\log{\Phi_{\beta}(N)}}}\right)
\end{eqnarray*}
Finally, on the whole integration space, we get the following upper bound :
\begin{eqnarray}\label{step3}
  {\EE}_{X}\sup_{(\lambda_{-J},\Lambda_{-J}),\cdots,(\lambda_{J},\Lambda_{J})\in H^{2J+1}}\sup_{-J\leq j\leq J}{\EE}_{g,g^{'}}
    \sup_{\alpha\in I_M}\vert G_{j}(\alpha)\vert\\
    \leq C\left(1+{\EE}_X\sup_{N\geq 1}\sqrt{\frac{\log^{+}{(\mid X_N\mid)}}
    {\log{\Phi_{\beta}(N)}}}\right)<\infty
\end{eqnarray}

\begin{figure}[h]
\psfrag{A}{$\lambda$}
\psfrag{B}{$\Lambda$}
\psfrag{h}{h}
\centerline{\includegraphics[scale=0.7]{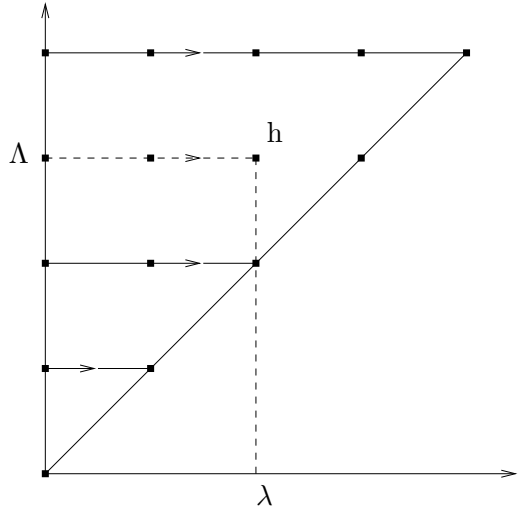}}
\caption{renumbering the upper triangle in ${\N}^2$}
\label{renumbfig}
\end{figure}

{\bf -Step 4 : } To end the proof of theorem \ref{principal}, it remains to deal with
  $$\sup_{-J\leq j \leq J}{\EE}_{g,g^{'}}\sup_{(\lambda,\Lambda)\in H}\sup_{\alpha\in I_M\cap {\mathbb Q}}G(j,\lambda,\Lambda,\alpha)$$
Let us fix $j$. In order to apply inequality \ref{ineg1}, we need to replace the supremum on $(\lambda,\Lambda)\in H$ by a supremum on only one
variable, in other words, we need to renumber $H$ using one variable $h\in{\N}^2$ given by the formula :
  $$h=\lambda+\frac{\Lambda(\Lambda-1)}{2}$$
(see figure \ref{renumbfig}).\\

Let us fix $j$. Inequality \ref{ineg1} gives us :
\begin{eqnarray}\label{renum}
  {\EE}_{g,g^{'}}\sup_{h\geq 1}\sup_{\alpha\in I_M\cap {\mathbb Q}}G(j,h,\alpha)\le K_1(\sup_{h\geq 1}{\EE}_{g,g^{'}}
  \sup_{\alpha\in I_M\cap {\mathbb Q}}G(j,h,\alpha)+{\EE}_{g,g^{'}}\sup_{h\geq 1}\vert\lambda_h\sigma_h\vert)
\end{eqnarray}
where $(\lambda_h)_{h\geq 1}$ is an isonormal sequence. The inequality :
  $$\frac{\Lambda(\Lambda-1)}{2}\le h\le\frac{\Lambda(\Lambda+1)}{2}$$
gives a polynomial dependence between $\Lambda$ and $h$, hence :
  $$\sigma_h={\mathcal O}(\frac{1}{\sqrt{\log h}})$$
and the first term in the right hand side of inequality (\ref{renum}) is dealt with in the same way as before. Finally, we get :
  $${\EE}_{g,g^{'}}\sup_{(\lambda,\Lambda)\in H}\sup_{\alpha\in I_M\cap{\mathbb Q}}G(j,\lambda,\Lambda,\alpha)
    \leq C\left(1+\sup_{N\geq 1}\sqrt{\frac{\log^{+}(\vert X_{N}\vert)}{\log{\Phi_{\beta}(N)}}}\right)$$
That is to say, by integrating on the whole space,
\begin{equation}\label{step4}
  {\EE}_{X}\sup_{1\leq j\leq n}{\EE}_{g,g^{'}}\sup_{(\lambda,\Lambda)\in H}\sup_{\alpha\in I_M\cap{\mathbb Q}}
    G(j,\lambda,\Lambda,\alpha)\leq C\left(1+{\EE}_{X}\sup_{N\geq 1}\sqrt{\frac{\log^{+}(\vert X_{N}\vert)}{\log{\Phi_{\beta}(N)}}}\right)
\end{equation}
Let us prove that 
  $${\EE}_X\sup_{N\geq 1}\sqrt{\frac{\log^{+}(\vert X_{N}\vert)}{\log{\Phi_{\beta}(N)}}}<\infty$$
Using Jensen inequality, we can get rid of the square root.
Let $\delta >0$. For all $N\geq 1,$ we have
  $$\beta\log ^{+}\vert X_{N}\vert \leq \beta\log ^{+}\left[\frac{\vert X_{N}
    \vert}{\Phi_{\beta}^{\delta}({N})}\right]+\beta\log ^{+} [\Phi_{\beta}^{\delta}({N})]$$
noticing that $\Phi_{\beta}^{\delta}({N}))\geq 2$ it is sufficient to show
  $${\EE} \sup_{N \geq 1} \log^{+}\left[\frac{\vert X_{N}\vert^\beta}{\Phi_{\beta}^{\delta\beta}(N)}
    \right]<\infty$$
Using now the inequality $\log^{+}(x)\leq x$ for any $x\geq 0$, it is sufficient to 
prove
  $$\sum_{N \geq 1}\frac{{\EE}\vert X_{N}\vert^{\beta}}{\Phi_{\beta}^{\delta\beta}(N)}<\infty$$
And as ${\EE} \vert X_N\vert^{\beta}\le\Phi_{\beta}(N)$ and $\Phi_{\beta}(N)\geq N$, if we chose $\delta=\frac{3}{\beta}$
we get the conclusion.
Steps  3 (see \ref{step3}), step 4 (see \ref{step4}) lead us to the announced result of theorem \ref{principal}.

\section{Applications}

Let us begin by giving an example where the $(X_k)$ are uniformly distributed :

\begin{ex}
Suppose that ${\mathcal L}(X_k)={\mathcal U}([\mu_k-\sigma_k/2,\mu_k+\sigma_k/2])$ with $\sigma_k>0$ et ${\EE}(X_k)=\mu_k$ with 
$\mu_k={\mathcal O}(k^d)$
for some $d>0$. The
characteristic function of $X_k$ can easily be computed :
  $$\varphi_{X_k}(t)=\frac{e^{2i\pi t\mu_k}}{\pi t\sigma_k}\sin(\pi t\sigma_k)$$
Using condition ${\mathcal H}$ of theorem \ref{moment}, the following condition 
\begin{equation}\label{condh}
  \sum_{n\geq 1}\frac{\vert a_n\vert}{\sigma_n}\mbox{\ converges}
\end{equation}
and 
  $$\sum_{n\geq 1}\frac{\sqrt{\sum_{k\geq n}\vert a_k\vert^2}}{n\sqrt{\log n}}<+\infty$$
are sufficient to get the desired convergence.\\
Notice that using corollary \ref{separationcauchy}, condition (\ref{condh}) is replaced by 
  $$\sum_{n\geq 1}\frac{1}{\sigma_n^2}<+\infty$$
The subexponential case could be dealt with in the same way.\\
If we consider the border case $a_k={\mathcal O}(k^{-1/2-\varepsilon})$, it is sufficient that :
  $$\exists\eta>0, \sigma_k\geq k^{\frac{1}{2}+\eta}$$
in this case :
  $${\PP}\{\forall k,X_k\in[\mu_k-\frac{k^{\frac{1}{2}+\eta}}{2},\mu_k+\frac{k^{\frac{1}{2}+\eta}}{2}]\}=1$$
which gives an information on the possible dispersion of the variables $X_k$.
\end{ex}

Here are other examples where the conditions of our theorems can be quite easily verified.    

\begin{cor}\label{trans}
Let $(X_k)_{k\geq 1}$ be a sequence of real independent random variables
whose law can be written in the following way for all $k\geq 1:$ ${\mathcal L}(X_k)=
{\mathcal L}(\sigma_k\cdot X+\mu_k)$ where $X$ verifies ${\EE}\vert X\vert^{\beta}<\infty$ for some $\beta>0$. Moreover,
we assume that there exist $d>0$ and $\delta>0$ such
that $\vert\sigma_k\vert={\mathcal O}(k^d)$, $\vert\mu_k\vert={\mathcal O}(k^d)$, the application
$t\mapsto t^{\delta}{\EE}\exp{(2i\pi tX)} $ is bounded on ${\R},$ and
let $(a_k)_{k\geq 1}$ be a sequence of real or complex numbers satisfying the following two conditions
\begin{enumerate}
\item $\vert a_k\vert=O(k^{-\beta})$ with $\beta>1/2$
\item $\sum_{k=1}^{\infty}\frac{1}{\vert \sigma_k\vert^{2\delta}}<\infty$
\end{enumerate}
Then there exists a measurable set $\Omega_o$ with full measure (${\PP}(\Omega_o) =1$)
such that for any $\omega \in\Omega_o$ for all $f\in B({\T})$ such that $\int_{{\T}}f(t)dt=0$ : 
for any compact $K$ which does not contain $0$,
the application $t\in K\mapsto F(t)=\sum_{k\geq 1}a_k f{(tX_k(\omega))}$  is
continuous and the series defining $F$ converges uniformly on $K$.
\end{cor}

The proof of corollary \ref{trans} relies on corollary \ref{separationcauchy}.

\begin{ex}The random variable $X$ may have a gaussian law with mean zero and variance one, 
a Cauchy law, the first Laplace law, an exponential law with parameter $\lambda>0.$ Let us precise the gaussian case.\\
Here ${\mathcal L}(X)={\mathcal N}(0,1)$, we have ${\EE}\exp{2i\pi tX}=e^{-t^2/2}$. Hence we can use the fact that
$t\mapsto e^{t^2/2}{\EE}\exp{2i\pi tX}$ is bounded on ${\R}$. In this case, the sufficient condition to obtain convergence is
  $$\forall\varepsilon>0, \exists N>0, \sup_{m>n\geq N}\sup_{\alpha\in K}\sup_{j\in {\Z} -\{0\}}\vert 
  \sum_{k=n}^m a_ke^{2i\pi\alpha\mu_kj}e^{-j^2\alpha^2\sigma_k^2/2}\vert<\varepsilon$$
Let $d(0,K)$ be the distance between $0$ and the compact $K$. Using the fact that $\vert a_k\vert=O(k^{-\beta})$ with $\beta>1/2$, 
the previous condition will be satisfied as soon as :
  $$\exists\varepsilon>0, \sigma_k\geq \frac{\sqrt{2(1-\beta+\varepsilon)}}{d(0,K)}\sqrt{\log k}$$
which, in terms of dispersion of the variables $X_k$ means that infinitely often:
  $$X_k\in[\mu_k-3\frac{\sqrt{2(1-\beta+\varepsilon)}}{d(0,K)}\sqrt{\log k},
  \mu_k+3\frac{\sqrt{2(1-\beta+\varepsilon)}}{d(0,K)}\sqrt{\log k}]$$
\end{ex}

\par
We discuss now the case when the laws of $X_k$ are generated by a convolution product of a given law $\mu$.
We distinguish two cases : on one hand when the support of $\mu$ contains non integer values, on the other hand when 
the support of $\mu$ is contained in ${\Z}$.
The first case is discribed by the following corollary :

\begin{cor}\label{aperiodic}
Let $(X_k)_{k\geq 1}$ be an sequence of real valued independent random variables
such that for all integer $k\geq 1$, ${\mathcal L}{(X_k)}=\mu^{*k}$ where $\mu$ is a probability measure
on ${\R}$ with ${\EE}\vert X_1\vert^{\delta}<\infty$ for some $\delta$.
Assume the following :
  $$\varphi_{X_1}(t)=1\Longleftrightarrow t=0\ \ \ \ \mbox{($X_1$ aperiodic)}\leqno{(a)}$$
  $$\sup_{t\in{\R}}|t^{\delta}{\EE}\exp{(2i\pi tX_1)}|=q<\infty\leqno{(b)\quad\exists\delta>0,}$$
Let $(a_k)_{k\geq 1}$ be a sequence of real or complex numbers such that the sequence
$\vert a_k\vert$ is decreasing and fulfills the two following conditions :
\begin{enumerate}
\item $\vert a_k\vert=O (k^{-\beta})$ avec $\beta>1/2$
\item $\sum_{k=1}^{\infty}\vert a_k-a_{k+1}\vert<\infty$
\end{enumerate}
Then there exists a measurable set $\Omega_o$ with full measure (${\PP}(\Omega_o) =1$)
such that for any $\omega \in \Omega_o$, for all $f\in B({\T})$ such that $\int_{{\T}}f(t)dt = 0$,
for any compact $K$ which does not contain $0$,
the application  $t\in K\mapsto F(t,\omega)=\sum_{k\geq 1}a_k f{(t\, X_k(\omega))}$ is continuous and
the series defining $F$ converges uniformly on $K$.
\end{cor}

\begin{rem}
If $X_1$ is striclty aperiodic ($\vert\varphi_{X_1}(t)\vert=1\Longleftrightarrow t=0$), then the condition on the differences
$\vert a_k-a_{k+1}\vert$ may be removed, using corollary \ref{separationcauchy} and the same kind of conputation as in the following proof.
\end{rem}

The random variable $X_1$ being real valued, its characteristic function is not periodic.
Take for example a gaussian law with mean zero and variance one.

\begin{dem}
Let $K$ be a compact which does not contain $0$. Using Abel 's summation, it is sufficient to prove
  $$\sup_{t\in K}\sup_{|j|\geq 1}\left\vert\sum_{k=1}^N({\EE}\exp{2i\pi jtX_1})^k\right\vert<\infty$$
independently of $N$.
Let us split the supremum on $j$ respectively into the supremum on the
indexes $J(q)$ and $\bar{J}(q)$
where $J(q)=\{j\in{\Z}^{*}:|j|^{\delta}\leq\left[\frac{2q}{\varepsilon^{\delta}}\right ]\}$ and $2\varepsilon$
is the distance between $0$ and the fixed compact $K$.\\
On one hand, using $(a)$, it can be proved that :
  $$\inf_{\vert t\vert>\varepsilon}\vert t\vert^{\delta}\vert 1-{\EE}\exp(2i\pi tX_1)\vert>0\leqno{\forall\varepsilon>0}$$
this implies
  $$\sup_{t\in K}\sup_{j\in J(q)}\left\vert\sum_{k=1}^N({\EE}\exp{2i\pi jtX_1})^k\right\vert\leq
    \sup_{t\in K}\sup_{j\in J(q)}C(\varepsilon)\vert jt\vert^{\delta}\leq C(K)\left[\frac{2q}{\varepsilon^{\delta}}\right]$$
On the other hand, using $(b)$, 
\begin{eqnarray*}
  \sup_{t\in K}\sup_{j\in \bar{J}(q)}\left\vert\sum_{k=1}^N({\EE}\exp{2i\pi jtX_1})^k\right\vert&\leq&
    \sup_{t\in K}\sup_{j\in\bar{J}(q)}\sum_{k=1}^N\left(\frac{q}{\vert tj\vert^{\delta}}\right)^k\\
  &\leq& C\sum_{k=1}^N\frac{1}{2^k}\leq 2C
\end{eqnarray*}
where $C$ is a universal constant.\\
\end{dem}

As for the integer valued case, we have :

\begin{cor}\label{integer}
Let $(X_k)_{k\geq 1}$ be an sequence of integer valued independent random variables
such that for all integer $k\geq 1$, ${\mathcal L}{(X_k)}=\mu^{*k}$ where $\mu$ is a probability measure
on ${\R}$ with ${\EE}\vert X_1\vert^{\beta}<\infty$ for some $\beta>0$. Let $(a_k)_{k\geq 1}$ be a sequence of complex numbers
such that $\vert a_k\vert=O(k^{-\beta})$ with $\beta>1/2$.\\
Assume either :
  $$\varphi_{X_1}(t)=1\Longleftrightarrow t=0\mbox{\ \ \ and\ \ \ }\sum_{k\geq 1}\vert a_k-a_{k+1}\vert\mbox{\ \ converges}$$
or :
  $$\vert\varphi_{X_1}(t)\vert=1\Longleftrightarrow t=0$$
Then there exists a measurable set $\Omega_o$ with full measure (${\PP}(\Omega_o) =1$)
such that for any $\omega \in \Omega_o$, for all $f\in A({\T})$ such that $\int_{{\T}}f(t)dt = 0$,
for any compact $K$ of the torus which does not contain $0$,
the application  $t\in K\mapsto F(t,\omega)=\sum_{k\geq 1}a_k f{(t\, X_k(\omega))}$ is continuous and
the series defining $F$ converges uniformly on $K$.
\end{cor}

\begin{ex}
If the law of $X_1$ is a Poisson law with parameter $1$, we use :
  $$\forall t\in{\T}, \vert\varphi_{X_1}(t)\vert\le e^{cos(2\pi t)-1}$$
\end{ex}

\noindent Dominique Schneider \\
Universit\'e du Littoral C\^ote d'Opale \\
L.M.P.A. CNRS EA 2597\\
50, rue F.Buisson B.P. 699 
F-62228 Calais cedex\\
dominique.schneider@lmpa.univ-littoral.fr \\ \\
Fr\'ed\'eric Paccaut \\
Universit\'e de Picardie Jules Verne \\
L.A.M.F.A. CNRS UMR 6140 \\
33, rue Saint Leu \\
F-80039 Amiens cedex 01 \\
frederic.paccaut@u-picardie.fr

\end{document}